\documentclass[a4paper]{article}

\usepackage{amsmath}

\usepackage{amsthm}
\usepackage[dvips]{graphicx}

\usepackage{epsfig}

\usepackage{amssymb}
\usepackage{latexsym}

\newtheorem{dfn}{Definition}[section]
\newtheorem{cor}[dfn]{Corollary}
\newtheorem{ex}[dfn]{Example}
\newtheorem{lem}[dfn]{Lemma}
\newtheorem{prop}[dfn]{Proposition}
\newtheorem{rem}[dfn]{Remark}
\newtheorem{thm}[dfn]{Theorem}

\begin{document}

\title{Rigidity of representations in $SO(4,1)$\\
for Dehn fillings on 2-bridge knots}
\author{Stefano Francaviglia and Joan Porti}
\date{\today}

\maketitle

\begin{abstract}
We prove that, for a hyperbolic two-bridge knot, infinitely many Dehn fillings are rigid in $SO_0(4,1)$.
Here rigidity means that any  discrete and faithful representation in $SO_0(4,1)$ 
is conjugate to the holonomy representation in $SO_0(3,1)$. We also show local rigidity for almost all Dehn fillings.
\end{abstract}

\tableofcontents

\section{Introduction}

Along this paper we consider compact, orientable three-manifolds $M$, whose
boundary is a torus $\partial {M}\cong T^2$ and whose interior admits a
complete hyperbolic metric of finite volume. 
We will specialize to the case where $M$ is the exterior of a hyperbolic two-bridge knot,
in particular its fundamental group is generated by two peripheral elements.

We tacitly assume that a basis for $H_1(\partial M;\mathbf Z)\cong\mathbf Z^2$
has been fixed, and for $(p,q)\in \mathbf Z^2$ coprime, we denote by
$M_{p/q}$ the manifold obtained by Dehn filling with meridian
curve $(p,q)$. 
According to Thurston's hyperbolic Dehn filling theorem, for all but
finitely many
 $p/q\in\mathbf Q\cup\{\infty\}$, the Dehn filled manifold 
$M_{p/q }$ is hyperbolic. In particular, its holonomy representation is 
the only discrete and faithful representation of $\pi_1(M_{p/q})$ in
$SO_0(3,1)$ up to conjugation. Here $SO_0(3,1)$ denotes the identity
component of $SO(3,1)$ and is isomorphic to
$\operatorname{Isom}^+(\mathbf H^3)$,  
 the orientation preserving isometry group of hyperbolic three space.

\vskip\baselineskip
In this paper we address the question of whether  
 $M_{p/q}$ has other discrete and faithful representations 
in $SO_0(4,1)\cong \operatorname{Isom}^+(\mathbf H^4)$.

For the unfilled manifold $M$, M. Kapovich proved global
rigidity~\cite{KapovichArticle}. That is to say, every discrete
and faithful representation 
of $\pi_1(M)$ in $SO_0(4,1)$ is conjugated to a representation in
$SO_0(3,1)$, and therefore to the holonomy of the
hyperbolic metric. Moreover M. Kapovich \cite{KapovichArticle}
 proved infinitesimal rigidity for infinitely
many Dehn fillings $M_{p/q}$, which are therefore locally rigid: 
there is no continuous non-trivial deformation in $SO_0(4,1)$
of the holonomy of $M_{p/q}$. That result was then generalized by
K. Scannell to Dehn fillings on a larger class of manifolds
\cite{ScannellSurgery}.  

Here we prove global rigidity for infinitely many Dehn
fillings on $M$, and local rigidity for almost all of them.

\begin{thm} \label{thm:limitslope}
Let ${M}$ be the exterior of a hyperbolic two-bridge knot.
Then, for infinitely many  $p/q\in\mathbf Q\cup\{\infty\}$, the Dehn filled manifold
$M_{p/q}$ has no discrete and faithful representation in $SO_0(4,1)$ other
than its holonomy in $SO_0(3,1)$.
\end{thm}

The following definition can be found in \cite{ScannellSurgery}. Here we consider the action of the holonomy representation in the Lie algebra
$\mathfrak{so}(4,1)$ via the adjoint representation.

\begin{dfn}
Let $M$ be a compact three manifold with boundary consisting of  tori and  whose interior is hyperbolic. 
Its \emph{parabolic cohomology} is defined to be the kernel
$$
PH^1(M,\mathfrak{so}(4,1))=\ker ( H^1(M,\mathfrak{so}(4,1))\to H^1(\partial M,\mathfrak{so}(4,1))   )
$$
\end{dfn}

In \cite{ScannellSurgery} Scannell  proves that, for two-bridge knot exteriors, $PH^1(M,\mathfrak{so}(4,1))=0$.
The following improves a theorem of \cite{ScannellSurgery}:

\begin{thm} \label{thm:localrigidity}
Let ${M}$ be a compact three manifold with boundary a torus whose interior is hyperbolic. If 
$PH^1(M,\mathfrak{so}(4,1))=0$,
then almost all  Dehn filled manifolds
$M_{p/q}$ are locally rigid in $SO(4,1)$.
\end{thm}

 Theorem~\ref{thm:limitslope} follows from Theorems~\ref{t_c} and \ref{prop:limitslope} below. 
Theorem~\ref{thm:localrigidity} is a consequence of the local analysis of the variety of representations made in 
the proof of Theorem~\ref{prop:limitslope}.

Convergence of representations will be understood in the variety of representations of $M$,
 in particular algebraic convergence.
For a discrete group $\Gamma$,  the set of all representations of $\Gamma$ in $SO_0(n,1)$ is denoted by  
$
R(\Gamma, SO_0(n,1))
$. 
It is well known that it is a real algebraic variety, cf.~\cite{JohnsonMillson,Morgan}.
We are interested in representations up to conjugacy, but the space of conjugacy classes 
does not seem to have a structure easy to work with. For results about this  set of conjugacy classes 
we refer to \cite{JohnsonMillson, KapovichArticle,Morgan,Scannell}.

When $\Gamma=\pi_1(M)$, it is customary to write 
$$
     R(M,SO_0(n,1))=R( \pi_1(M),SO_0(n,1)).
$$

\begin{thm}\label{t_c}
  Let $M$ be  a  hyperbolic two-bridge knot exterior. Let
  $\rho_0:\pi_1(M)\to{\rm Isom}(\mathbf H^4)$ be the holonomy of the
  (Fuchsian) complete hyperbolic structure of $M$. 
Let $\{\rho_n\}_{n\in\mathbf N}$ be
  a sequence of representations
$$\rho_n:\pi_1(M)\to{\rm Isom}(\mathbf H^4).$$ Suppose that, for
each $n$, the representation $\rho_n$
factorizes through a discrete and faithful representation of $\pi_1(M_{p_n/q_n})$
with $p_n^2+q_n^2\to\infty$.
Then, up to conjugation, $\rho_n$ converges to $\rho_0$.
\end{thm}

\begin{thm} \label{prop:limitslope}
Let ${M}$ be hyperbolic with one cusp with $PH^1(M,\mathfrak{so}(4,1))=0$.
Let $\rho_n$ be a sequence of representations of $\pi_1(M)$ in $SO_0(4,1)$ 
that converges to the holonomy of $M$.

If each $\rho_n$ comes from a discrete and faithful representation of
$M_{p_n/q_n}$, not conjugated to a representation in $SO_0(3,1)$, then 
$p_n/q_n\to l$ for some $l\in \mathbf R\cup\{\infty\}$ depending only on M.
\end{thm}

Here is the idea of the proof of  Theorem~\ref{t_c}.
Take a sequence $\{\rho_n\}_{n\in\mathbf N}$, so that each $\rho_n$ is a discrete and faithful representation
of $M_{p_n/q_n}$ with $p_n^2+q_n^2\to \infty$. By theorems of Bestvina,  Morgan-Shalen and Paulin, 
the space of discrete and faithful representations of $M$
in $SO(4,1)$ is compact  \cite{Bestvina,MorganShalen, Paulin}. 
Even if $\rho_n$ is not a faithful representation of $M$,
the proof can be adapted to say that a
subsequence of $\rho_n$ converges to $\rho$, a representation of $M$.
Moreover, the so called Chuckrow-Wielenberg's theorem can also be adapted to say that
$\rho$ is discrete and faithful. 
Since $M$ is generated by two peripheral elements, by a theorem of M.~Kapovich \cite{KapovichArticle} $\rho$ must be the holonomy representation of $M$ in 
$SO(3,1)$.

Theorem~\ref{prop:limitslope} involves an analysis of $R(M,SO(4,1))$ in a neighborhood of $\rho_0$,
following closely the results obtained by Scannell in \cite{ScannellSurgery}.
The tangent space to 
$$
R(\pi_1 M, SO_0(4,1))
$$ 
is the space of cocycles valued on the Lie algebra $\mathfrak{so}(4,1)$, 
and the space tangent to the orbits by conjugation is the space of coboundaries. Thus, to get relevant information, we study
 the cohomology group 
$H^1(M,\mathfrak{so}(4,1))$. We need to understand how a sequence of representations in $SO(4,1)$ can approach $\rho_0$, and we shall 
study which elements in $H^1(M,\mathfrak{so}(4,1))$ are tangent vectors of deformations of $\rho_0$ in $R(\pi_1(M), SO_0(4,1))$, what in Subsection~\ref{subsectio:differentialtgcone} is called the
\emph{differentiable tangent cone}.

 Since $\rho_0$ is contained in $SO(3,1)$, the Lie algebra splits as
 $\pi_1(M)$-module as follows
$$\mathfrak{so}(4,1)=\mathfrak{so}(3,1)\oplus\mathbf R^{3,1},
$$
where $\mathbf R^{3,1}$ denotes the Minkowski space equipped with the linear
action of $SO(3,1)$.
That splitting induces a direct sum of cohomology groups
$$
H^1(M,\mathfrak{so}(4,1))=H^1(M,\mathfrak{so}(3,1))\oplus H^1(M,\mathbf R^{3,1}).
$$
The subspace $H^1(M,\mathfrak{so}(3,1))$ has dimension 2 and it is the tangent space to the variety of representations in 
$SO(3,1)$ up to conjugation, described by Thurston's proof of the hyperbolic Dehn filling theorem, cf.~\cite{KapovichBook}. 
By a theorem of Scannell \cite{ScannellSurgery}, $H^1(M,\mathbf R^{3,1})$ has dimension one. In Proposition~\ref{prop:tangentcone},
we show that the tangent cone is contained in the union
$$
H^1(M,\mathfrak{so}(4,1))=H^1(M,\mathfrak{so}(3,1))\cup  H^1(M,\mathbf R^{3,1}).
$$
With the help of the curve selection lemma, this implies that if we
have a sequence of non-Fuchsian representations 
approaching $\rho_0$, then the sequence  must be contained in a semialgebraic set tangent to $H^1(M,\mathbf R^{3,1})$.

In the proof we need to understand how elements in $H^1(\partial M,\mathbf R^{3,1})$  are realized by deformations in the boundary $\partial M$,
viewed as isometries of $\mathbf R^3=\partial \mathbf H^4\setminus\{\infty\}$.
Those are realized by deforming a lattice in $\mathbf R^2\subset \mathbf R^3$ as a group of screw motions whose axis is contained in 
$\mathbf R^2$. This kind of deformation imposes some restriction on the Dehn filling coefficients, that we use to prove the theorem.

Our study of the variety of representations relies on previous work of Scannell~\cite{ScannellSurgery}, where 
he shows that this is a singular point of the variety of representations.

The paper is organized in four sections. In
Section~\ref{section:convergence} we prove
Theorem~\ref{t_c}. Section~\ref{sec:infinitesimal} is devoted to the
preliminaries about infinitesimal deformations,  and
Section~\ref{sec:localanalysis} to the analysis of a neighborhood of
the variety of representations. 
The results of both sections are used
 in Section~\ref{sec:sequencesreps}, where
 Theorem~\ref{prop:limitslope} is proved.

\paragraph*{Acknowledgments} We thank Aris Daniilidis and Francesco
Bonsante for useful conversations. The first author was partially
supported by the European Research Council (MERG-CT-2007-046557)
 and the second one by grant FEDER/MEC MTM2006-04353. The first author
 wishes to thank the Departament de Matem\`atiques de la UAB for the kind
 hospitality during his frequent visits.

\section{Convergence of representations}
\label{section:convergence}

In this section we prove Theorem~\ref{t_c}. Even if most of the arguments and techniques we use
 can be found
in the literature, we give a proof for completeness, stressing the changes required in 
our situation.
We mention in particular the recent work of
M. Kapovich~\cite{KapovichPreprint07082671} 
on convergence of groups.

\proof[Proof of Theorem~\ref{t_c}] The proof goes through two main steps:

\begin{itemize}
\item[--] \emph{Step 1.} The sequence $\{\rho_n\}$ is bounded in 
$R(\pi_1(M),{\rm Isom}(\mathbf H^4))$.
\item[--] \emph{Step 2.} Each accumulation point of $\{\rho_n\}$ is discrete and
  faithful, whence conjugate to $\rho_0$ by \cite{KapovichArticle, ScannellSurgery}.
\end{itemize}

Each step is explained in a different subsection.

\subsection{Ultralimits and asymptotic cones}

We concentrate here on the first step of the proof of
Theorem~\ref{t_c}. We start by recalling the
following result of Morgan and Shalen \cite{MorganShalenIII}.

\begin{thm}\label{t_ms}
  If $M$ is a complete hyperbolic $3$-manifold, then there is no
  action of $\pi_1(M)$ on an $\mathbf R$-tree by isometries having
  \begin{itemize}
  \item no global fixed points, and
  \item small arc-stabilizers ({\em i.e.} arc-stabilizers do not
    contain rank-two free sub-groups).
  \end{itemize}
\end{thm}

We will show that, if $\{\rho_n\}$ was unbounded, then it would induce an action
of $\pi_1(M)$ on an $\mathbf R$-tree as the one
forbidden by Theorem~\ref{t_ms}, deducing therefore that
$\{\rho_n\}$ must be bounded. The way to do that uses standard
techniques of asymptotic cones (see for example~\cite{KapovichBook}
for details on 
asymptotic cones.)

\ 

Let $\omega$ be a non-principal ultrafilter, that we think of as a family of
subsets of natural numbers such that:

\begin{enumerate}
\item Given any subset $S\subset \mathbf N$, either $S\in \omega$ or
  $\mathbf N\setminus S \in \omega$;
\item if $S\in \omega$, and $S^\prime \supset S$, then $S^\prime \in \omega$;
\item if $S \subset \mathbf N$ is a finite subset, then $S\notin \omega$;
\item if $S,S^\prime \in \omega$, then $S\cap S^\prime \in \omega$.
\end{enumerate}

We say that a sequence $\{x_n\}$ in a topological space $\omega$-converges to a
point $x$, and we write $\omega\lim x_n=x$, if for each open neighborhood $U$
of $x$, the set $\{i\in\mathbf N: x_i\in U\}$ belongs to $\omega$. It
is an easy exercise that any sequence in a compact space has a unique
$\omega$-limit.

\medskip

Let
$\{\gamma_i\}$ be a finite set of generators of $\pi_1(M)$. Let
$*=\{*_n\}$ be a sequence of points such that $*_n\in\mathbf H^4$ realizes
the minimum 
$$\min_{p\in\mathbf H^4}\max_{i}\,d(\rho_n(\gamma_i)p,p)$$
and let
$$\lambda_n=\max_{i}\,d(\rho_n(\gamma_i)*_n,*_n).$$

Let $(X,d)$ be the asymptotic cone of $\mathbf H^4$ made using the
ultrafilter $\omega$, the  sequence of rescaling parameters $\{\lambda_n\}$ and the
base-points sequence $\{*_n\}$. Namely, 
$$
X=\left\{\{x_n\}\subset \mathbf H^4 : \omega\lim\,
  \frac{d(x_n,*_n)}{\lambda_n}<\infty\right\}
$$ 
where we identify two sequences $\{x_n\}$ and $\{y_n\}$ whenever 
$$
\omega\lim\, \frac{d(x_n,y_n)}{\lambda_n}=0,
$$
and we set $$d(\{x_n\},\{y_n\})=\omega\lim\, \frac{d(x_n,y_n)}{\lambda_n}.$$

The following lemma is a standard fact about asymptotic cones of
hyperbolic spaces, see for example~\cite{KapovichBook}.

\begin{lem}
  If $\{\rho_n\}$ is unbounded -- namely, if $\lambda_n\to\infty$ --
  then $(X,d)$ is an $\mathbf R$-tree.
\end{lem}

\begin{lem}
  The $\omega$-limit $\rho_\omega$ of $\{\rho_n\}$ is an isometric action of
  $\pi_1(M)$ on $X$. That is to say, a representation of $\pi_1(M)$ on
  ${\rm Isom}(X)$.
\end{lem}

\begin{proof}
  The action $\rho_\omega$ on $X$ is tautologically defined by
$$\rho_\omega(\gamma)(\{x_n\})=\{\rho_n(\gamma)(x_n)\}.$$
We have to check that such a definition is well-posed; namely, that
for all $\{x_n\}\in X$ and all $\gamma\in\pi_1(M)$ we have
$\{\rho_n(\gamma)(x_n)\}\in X$. In other words, we need to check that 
$$\omega\lim\frac{d(x_n,*_n)}{\lambda_n}<\infty
\qquad\Longrightarrow\qquad
\omega\lim\frac{d(\rho_n(\gamma)x_n,*_n)}{\lambda_n}<\infty.$$

Let $\gamma=\gamma_{i_1}\cdots\gamma_{i_k}$ be a decomposition of
$\gamma$ in terms of the fixed generators of $\pi_1(M)$.
We have:
\begin{eqnarray*}
&&d(\rho_n(\gamma)x_n,*_n)\leq d(\rho_n(\gamma)x_n,\rho_n(\gamma)*_n)
+d(\rho_n(\gamma)*_n,*_n)\\
&=&d(x_n,*_n)+d(\rho_n(\gamma)*_n,*_n)\\
&\leq& d(x_n,*_n)+
\sum_{j=1}^k d(\rho_n(\gamma_{i_1}\cdots\gamma_{i_j})*_n,
\rho_n(\gamma_{i_1}\cdots\gamma_{i_{j-1}})*_n))\\
&\leq& d(x_n,*_n)+ k\max_i d(\rho_n(\gamma_i)*_n,*_n)= d(x_n,*_n)+k\lambda_n. 
\end{eqnarray*}

Therefore $\rho_\omega$ is well-defined. The fact that it is an action
by isometries is obvious.
\end{proof}

\begin{lem}
The representation $\rho_\omega\!:\pi_1(M)\to{\rm Isom}(X)$ 
has no global fixed point.
\end{lem}
\begin{proof}
Let $\{x_n\}$ be any point of $X$.
By construction of $*_n$ and $\lambda_n$, the index $i_n$ that realizes
realizes $\max_i\,d(\rho_n(\gamma_i)(x_n),x_n),$ satisfies:
$$
d(\rho_n(\gamma_{i_n})(x_n),x_n)\geq 
\lambda_n.$$
Thus
$$\frac{d(\rho_n(\gamma_{i_n})(x_n),x_n)}{\lambda_n}\geq1$$
and therefore  
$d(\rho_\omega(\gamma_\omega)\{x_n\},\{x_n\})\geq1$, where
$\gamma_\omega=\omega\lim\gamma_{i_n}\in\{\gamma_i\}$. 
It follows that $\{x_n\}\in X$ is
not globally fixed.
\end{proof}

Notice that the previous two lemmas, as well as the following one, 
do not use the hypothesis that
$\{\rho_n\}$ is unbounded. Indeed, such hypothesis is only needed to
show that $X$ is an $\mathbf R$-tree.

\begin{lem}
 The representation $\rho_\omega$
 has small arc-stabilizers. 
\end{lem}

\begin{proof}
Let $I\subset X$ be an arc, and let $\Gamma<\pi_1(M)$ be its
$\rho_\omega$-stabilizer: 
$$\rho_\omega(\Gamma)(I)=I.$$ 
Let $\{x_n\}$ and $\{y_n\}$ be the end-points of $I$, and let
$\gamma\in\Gamma$. Up to replacing  $\gamma$ by $\gamma^2$ we
have
$\rho_\omega(\gamma)\{x_n\}=\{x_n\}$ and $\rho_\omega(\gamma)\{y_n\}=\{y_n\}$ as
elements of $X$. That is to say,
$$
\omega\lim\frac{d(\rho_n(\gamma)(x_n),x_n)}{\lambda_n}=0
$$ 
and the same holds true for $y_n$. Since the $\omega$-limit of
$d(x_n,y_n)/\lambda_n$ is the length of $I$, there exists a
subsequence of indices $\{n_k\}$ such that 
$$\frac{d(\rho_{n_k}(\gamma)(x_{n_k}),x_{n_k})+d(\rho_{n_k}(\gamma)(y_{n_k}),y_{n_k})
}{d(x_{n_k},y_{n_k})}\to 0.
$$
Given any other $\psi\in\Gamma$, the same limit hods true
up to subsequence. Now we use Proposition~4.5
of~\cite{Bestvina} 
(a Margulis-type argument, together with some hyperbolic trigonometry) 
to deduce that the group generated by $\rho_{n_k}(\gamma)$ and
$\rho_{n_k}(\psi)$ is abelian for large enough $k$. 

It follows that the commutator $[\gamma,\psi]$ belongs to
the kernel of  $\rho_{n_k}$ for large enough $k$.

Since $\rho_n$ factorizes through a faithful
representation of $\pi_1({M_{p_n/q_n}})$, 
the following lemma shows that $[\gamma,\psi]$ is in fact trivial in 
$\Gamma$. This
implies that $\Gamma$ is abelian, and therefore cannot contain
rank-two free sub-groups.
\end{proof}

\begin{lem}\label{l_ker} Let $P_n\! :\! \pi_1(M)\to\pi_1({M_{p_n/q_n}})$ denote the natural surjection 
induced by Dehn filling. If $p_n^2+q_n^2\to\infty$, then, for all $m\in \mathbf N$, 
$$\bigcap_{n>m}\ker(P_n)=1.$$ 
\end{lem}

\begin{proof}
For large enough $n$, the Dehn-filling on $M$ with parameters
$(p_n,q_n)$ is hyperbolic by Thurston's hyperbolic Dehn-filling
theorem.
Let $\gamma\in\cap\ker(P_n)$. Then, the holonomy of $\gamma$ is trivial
in each ${M_{p_n/q_n}}$. Again, by Thurston's theorem, the holonomy of $\gamma$
in ${M_{p_n/q_n}}$ converges to the holonomy of $\gamma$ in  the complete hyperbolic structure of $M$, which is
therefore trivial. This is possible only if $\gamma=1$.
\end{proof}

So, $\rho_\omega$ is an isometric action on $X$ with no global fixed
points and small arc stabilizers. By Theorem~\ref{t_ms} such an action
cannot exist. Therefore $X$ cannot be an $\mathbf R$-tree, whence we
get that $\lambda_n$ must be bounded. Then, up to
conjugation by isometries of $\mathbf H^4$, we can suppose that $*_n$
is constant, and the sequence $\{\rho_n\}$ is in that
case bounded. This ends the proof of Step 1.

\subsection{Accumulation point of Dehn-fillings}

We deal now with the last step of the proof of Theorem~\ref{t_c}. 
Let $\rho$ be an accumulation point of
$\{\rho_n\}$. Let $V$ be an open neighborhood of the identity in ${\rm
Isom}(\mathbf H^4)$ such that any discrete group finitely generated by
elements in $V$ is virtually nilpotent. Such a $V$ exists by the Margulis
lemma.
Let $U$ be an open neighborhood of the identity such that $\overline
U\subset V$.

\begin{lem}\label{l_ab}
  Let $\Gamma<\pi_1(M)$ be the subgroup generated by the elements $\gamma$
  such that $\rho(\gamma)\in U$. Then $\Gamma$ is abelian.
\end{lem}
\begin{proof}
Let $\Gamma_0<\Gamma$ be a group finitely generated by elements whose
$\rho$-image is in $U$.
For $n$ large enough, $\rho_n(\Gamma_0)$ is a discrete subgroup of ${\rm
  Isom}(\mathbf H^4)$, finitely generated by elements of $V$. Then,
$\rho_n(\Gamma_0)$ is virtually nilpotent.

Since $\rho(\Gamma_0)$ is virtually nilpotent, it is also elementary, because 
 the limit set of a nontrivial normal subgroup is the same as the limit set of the whole group.
 In particular, $\rho_n(\Gamma_0)$ is elementary. Moreover, $\rho_n(\Gamma_0)$ is torsion-free 
because $\rho_n$ is faithful as a representation of $\pi_1(M_{p_n/q_n})$. 

An elementary and torsion-free group of isometries in $\mathbf H^4$ must be one of the following:
\begin{itemize}
\item either a subgroup of the stabilizer of a geodesic, $\mathbf R\rtimes O(3)$;
\item or a parabolic subgroup fixing a point at $\partial \mathbf H^4$, i.e.\ a subgroup 
of $\operatorname{Isom}(\mathbf R^3)\cong \mathbf R^3\rtimes O(3)$.
\end{itemize}

In particular, all such groups are nilpotent of order two, cf.~\cite{Wolf}. Thus, for $n$ large enough, 
$$\left[[\rho_n(\Gamma_0),\rho_n(\Gamma_0)],\rho_n(\Gamma_0)\right]=1.$$
It follows that any $\gamma\in\Gamma_0$ of the form
$\gamma=[[\gamma_1,\gamma_2],\gamma_3]$ belongs eventually to $\ker \rho_n$, and
by Lemma~\ref{l_ker}, this forces $\gamma$ to be trivial. That is to
say, $\Gamma_0$ itself is virtually nilpotent. Since $\Gamma_0$ is a
subgroup of the fundamental group of a hyperbolic manifold, this
implies that $\Gamma_0$ is abelian. Since this holds for any
$\Gamma_0$, we get that $\Gamma$ itself is abelian.
\end{proof}

\begin{cor}
  The representation $\rho$ is faithful.
\end{cor}

\begin{proof}
  If  $\ker(\rho)$ were not trivial, then it would be abelian by Lemma~\ref{l_ab}, 
but
$\pi_1(M)$ has no abelian, non-trivial, normal subgroups.
\end{proof}

The very same argument shows the following corollary.
\begin{cor}
  The representation $\rho$ is discrete.
\end{cor}

\begin{proof}
Let $H_0$ be the connected component of the identity of the 
topological closure $\overline{\rho(\pi_1(M))}$. 
As $\overline{\rho(\pi_1(M))}$ is a Lie
group, $H_0$ is normal in $\overline{\rho(\pi_1(M))}$. It follows that
$\Gamma_0:=H_0\cap\rho(\pi_1(M))$ is normal in $\rho(\pi_1(M))$, and since
$\rho$ is faithful, $\Gamma:=\rho^{-1}(\Gamma_0)$ is normal in
$\pi_1(M)$. 

On the other hand, $\Gamma_0$ is generated by elements in $U$, and
then, by Lemma~\ref{l_ab}, $\Gamma$ is abelian. Therefore,
$\Gamma=\{1\}$ and $H_0=\{1\}$, that is to say, $\rho$ is discrete.
\end{proof}

This concludes Step 2 and so the proof of Theorem~\ref{t_c}.\qed

\section{Infinitesimal isometries and deformations}
\label{sec:infinitesimal}

This section contains the background material and tools that we
need in the proof of Theorem~\ref{prop:limitslope}.

Let $\mathbf R^{n,1}$ denote the Minkowski space, i.e.\  $\mathbf R^{n+1}$ equipped
 with the usual Lorentz product, that has matrix
$$
J=\begin{pmatrix}
      -1 & 0 &  & \\
       0 & 1 &  &  \\
          &  & \ddots & \\
         &   &  &  1
  \end{pmatrix}.
$$
We shall use the hyperboloid model for hyperbolic space
$$
\mathbf H^n=\{x=\begin{pmatrix} x^0 \\ \vdots \\ x^n \end{pmatrix}\in \mathbf R^{n,1} \mid x^t J x=-1, x^0> 0\}
$$
so that  the orientation preserving isometry group of $\mathbf H^n$ is identified with the identity component
of the group of linear transformations of $\mathbf R^{n+1}$ that preserve $J$ as a bilinear form:
$$
\operatorname{Isom}^+(\mathbf H^n)=SO_0(n,1).
$$

We are interested in the cases $n=3$ and $n=4$. We shall consider the inclusion 
$SO_0(3,1)\subset SO_0(4,1)$ induced by the inclusion  $\mathbf R^{3,1}\subset \mathbf R^{4,1}$
consisting in adding a fifth coordinate $x^4$. Namely:
\begin{equation}
\label{eqn:canonicalinclusion} 
\begin{array}{rcl}
SO(3,1) & \to & SO(4,1) \\
A & \mapsto & \begin{pmatrix}
               A & 0 \\ 0 & 1
              \end{pmatrix}.
\end{array}
\end{equation}

\subsection{Infinitesimal isometries}

The Lie algebra of $SO(n,1)$ is
$$\mathfrak{so}(n,1)=\{a\in M_{n+1}(\mathbf R)\mid a^tJ=-J a\}.$$
 Elements in 
$\mathfrak{so}(n,1)$ are viewed as \emph{infinitesimal isometries}: a matrix $a\in \mathfrak{so}(n,1)$ is the tangent vector to the path
$\exp(t\, a)\in SO_0(n,1)$ at $t=0$. The action of the isometry group on itself by conjugation induces the adjoint action of 
$SO_0(n,1)$ on the Lie algebra $\mathfrak{so}(n,1)$. Since we have an inclusion $SO_0(3,1)\subset SO_0(4,1)$, $\mathfrak{so}(4,1)$ is also
a $SO_0(3,1)$-module.

 \begin{lem} 
\label{lem:descomposition}
 We have an isomorphism of $SO_0(3,1)$-modules
$$\mathfrak{so}(4,1)= \mathfrak{so}(3,1) \oplus \mathbf R^{3,1}
$$
where $SO(3,1)$ acts on $\mathfrak{so}(3,1)$ by the adjoint action and on $\mathbf R^{3,1}$ by the usual linear action.
 \end{lem}

\begin{proof} Explicit construction. Given a matrix $a\in \mathfrak{so}(3,1)$ and 
a (column) vector $v\in\mathbf R^{3,1}$, we consider the following matrix:
\begin{equation}
\label{eqn:r41inso41}
\begin{pmatrix}
a  & v \\
-v^t J & 0
\end{pmatrix} \in \mathfrak{so}(4,1),
\qquad 
\end{equation}
where $v^t$ is the transpose matrix.
It is easy to check that this gives the isomorphism of the lemma, compatible with inclusion (\ref{eqn:canonicalinclusion}).
\end{proof}

The Lie algebra of infinitesimal deformations can be identified with the Lie algebra of 
Killing vector fields.

\begin{prop}
\label{prop:R31orthogonal}
The subspace $\mathbf R^{3,1}\subset \mathfrak{so}(4,1)$ corresponds to the Killing vector fields orthogonal to $\mathbf H^3\subset \mathbf H^{4}$.
\end{prop}

\begin{proof}
Given $a\in \mathfrak{so}(4,1)$, the corresponding field $V$ evaluated at $x\in \mathbf H^4$ is 
$$
V_x=\left.\frac{d\phantom{t}}{dt}\right\vert_{t=0}\exp(t\, a) x.
$$
Since we are working in a linear model,
so that $x\in \mathbf H^4\subset\mathbf R^{4,1}$:
$$
V_x= a\, x\in \mathbf R^{4,1}.
$$
In this model, $\mathbf H^ 3=\mathbf H^4\cap \{ x^4=0\}$.
Thus the Killing vector fields perpendicular to  $\mathbf H^ 3$ correspond to matrices in 
$\mathfrak{so}(4,1)$ whose entries are zero, except for the last column or the last row, which is the image of the embedding of
$ \mathbf R^{3,1} $ in $\mathfrak{so}(4,1)$.
\end{proof}

The splitting of Lemma~\ref{lem:descomposition} can be also understood by using the action on the \emph{de Sitter space}
$$
\mathbf S^{3,1}=\{v\in \mathbf R^{4,1}: v^tJv=1\},
$$
which is naturally identified to the space of oriented hyperplanes in $\mathbf H^4$, cf.\ \cite{EpsteinPenner}.
Since $SO(4,1)$ acts transitively on $\mathbf S^{3,1}$ with
stabilizers $SO(3,1)$, we have that
$$SO(4,1)/SO(3,1)=\mathbf  S^{3,1}.$$
Moreover the fibration
$$
SO(3,1)\to SO(4,1)\to\mathbf S^{3,1}
$$
whose projection maps $A\in SO(4,1)$ to $A\cdot p$ for some fixed $p\in\mathbf S^{3,1}$, 
induces an exact sequence of $\mathfrak{so}(3,1)$-modules:
$$
0\to \mathfrak{so}(3,1)\to \mathfrak{so}(4,1)\to T_ p\mathbf S^{3,1}\to 0.
$$
which splits, by using either the Killing form or Lemma~\ref{lem:descomposition}.

\begin{rem}
 We have a canonical identification
$$
T_ p\mathbf S^{3,1}\cong \mathbf R^{3,1}
$$
where $p\in \mathbf S^{3,1}$ is the hyperplane stabilized by $\mathfrak{so}(3,1)\subset  \mathfrak{so}(4,1)$.
\end{rem}

In order to be compatible with the previous computations,  we shall assume that
we choose the point
$$
p=\begin{pmatrix} 
0\\ 0 \\ 0 \\ 0\\ 1   
  \end{pmatrix} \in\mathbf S^{3,1}.
$$

\subsection{The Zariski tangent space to the variety of representations}

For a representation $\rho\!:\Gamma\to SO_0(n,1)$, the Zariski
tangent space of the variety of representations
at $\rho$ is naturally identified to the space of cocycles \cite{Weil}.
Namely, the space of cocycles is defined as
$$
Z^1(\Gamma,\mathfrak{so}(n,1)_{\rho})=\{d:\Gamma\to \mathfrak{so}(n,1)\mid
d(\gamma_1\gamma_2)=d(\gamma_1)+Ad_{\rho(\gamma_1)}(d(\gamma_2))\}.
$$
\emph{Weil's correspondence} maps the cocycle $d\in Z^1(\Gamma,\mathfrak{so}(n,1)_{\rho})$ 
to the infinitesimal deformation
$$
\rho_t(\gamma)=(1+t\, d(\gamma))\rho(\gamma)\qquad \forall\gamma\in\Gamma,
$$ 
which is a representation mod $t^ 2$, hence  a Zariski tangent vector, cf.~\cite{LubotzkyMagid}.

The space of coboundaries is 
\begin{multline*}
B^1(\Gamma,\mathfrak{so}(n,1)_{\rho})=\{d_a:\Gamma\to \mathfrak{so}(n,1)\mid
d_a(\gamma)= a-Ad_{\rho(\gamma)}(a) \\ \textrm{ for some }a\in \mathfrak{so}(n,1)\}
\end{multline*}
and it is identified to the Zariski tangent space of orbits by conjugation.

The quotient is the cohomology group:
$$
H^1(\Gamma,\mathfrak{so}(n,1)_{\rho})=Z^1(\Gamma,\mathfrak{so}(n,1)_{\rho})/B^1(\Gamma,\mathfrak{so}(n,1)_{\rho}).
$$

Under some circumstances,  $H^1(\Gamma,\mathfrak{so}(n,1))$  can be viewed as the tangent space to the space of 
conjugacy classes of representations. However,  for technical reasons, since $\mathbf R$ 
is not algebraically closed, it is easier
to work with the variety of representations.


Now we focus on the case $\Gamma=\pi_1(M)$, where $M$ is a cusped manifold and $\rho_0$ is the holonomy representation of its complete structure.
We omit 
the representation $\rho_0$ when writing the Lie algebra as $\pi_1(M)$-modules via the adjoint action of $\rho_0$. We also write $H^1(M,V)$ for $H^1(\pi_1(M),V)$.
Recall that the parabolic cohomology is defined as:
$$
PH^1(M,\mathfrak{so}(4,1) )=\ker (    H^ 1(M,\mathfrak{so}(4,1) )\to H^1(\partial M, \mathfrak{so}(4,1) )     ).
$$
Next lemma is due to Scannell \cite{ScannellSurgery}. 
\begin{lem}\label{lem:kerboundary} 
Let $\rho_0:\pi_1(M)\to SO_0(3,1)\subset SO_0(4,1)$ be the
holonomy representation of a hyperbolic two-bridge knot. 
 Then 
$$
 PH^ 1(M,\mathfrak{so}(4,1) )    =0   .
$$
\end{lem}

To prove that $\ker ( H^ 1(M,\mathfrak{so}(4,1) )\to H^1(\partial M, \mathfrak{so}(4,1) )   )$ vanishes,
the idea is that an element in the kernel corresponds to a deformation that keeps the generators parabolic.
By a geometric argument due to Kapovich and Scannell, 
such a representation must preserve a hyperbolic space of dimension three.

\begin{lem}
\label{cor:halfdimension} Suppose that $PH^1(M,\mathfrak{so}(4,1))$ vanishes. 
Then the image of the inclusion
$$
0\to H^ 1(M,\mathfrak{so}(4,1)) \to H^ 1(\partial M,\mathfrak{so}(4,1))
$$
has half dimension.
\end{lem}

The proof consists in applying the long exact sequence of the pair and Poincar\'e duality.

Recall from Lemma~\ref{lem:descomposition} that,
as $\rho_0$-module by the adjoint action, we have a decomposition
$$
\mathfrak{so}(4,1)=\mathfrak{so}(3,1)\oplus \mathbf R^{3,1},
$$
where $\mathbf R^{3,1}$ is the four dimensional real vector space equipped with the linear action of
$SO(3,1)\subset GL(\mathbf R,4)$.
In particular:
$$
H^ 1(M,\mathfrak{so}(4,1) )=H^ 1(M,\mathfrak{so}(3,1) )\oplus H^ 1(M,\mathbf R^{3,1}).
$$

The dimensions of  those spaces for the torus are well known.

\begin{lem}[\cite{KapovichArticle,ScannellSurgery}]
\label{lem:dimH1}
$\dim (H^ 1(\partial M,\mathfrak{so}(3,1))) = 4$ and $\dim (H^ 1(\partial M,\mathbf R^{3,1})) = 2$.
\end{lem}

From Lemmas~\ref{lem:dimH1} and~\ref{cor:halfdimension}, we get:

\begin{cor} If $PH^1(M,\mathfrak{so}(4,1))=0$, then 
$\dim (H^ 1(M,\mathfrak{so}(3,1) )) = 2$ and $\dim (H^ 1(M,\mathbf R^{3,1})) = 1$.
\end{cor}

The following lemma can be easily proved using the formalism of the previous section. Recall that
the projection $SO(4,1)\to\mathbf S^{3,1}$ maps $A\in SO(4,1)$ to $A\cdot p\in \mathbf S^{3,1}$,
for some fixed $p\in\mathbf S^{3,1}$.

\begin{lem}\label{cor:cyclecbar}
  Let $\rho_t:\Gamma\to SO(4,1)$ be a smooth
  path of representations of a group $\Gamma$ in $SO(4,1)$, and let
  $\bar \rho_t$ be its projection to $\mathbf S^{3,1}$. 
Then, the first non-trivial derivative of $\bar \rho_t$ defines a
cocycle in $Z^1(\Gamma,\mathbf R^{3,1}).$

If the cocycle is nontrivial in $H^1(\Gamma,\mathbf R^{3,1})$, then for small $t$ the representation $\rho_t$ does not fix 
any  hyperplane in $\mathbf H^4$.
\end{lem}

Notice that this lemma uses the natural interpretation of points in $\mathbf S^{3,1}$ as hyperplanes in $\mathbf H^4$.

\subsection{Trace functions}

For $\gamma\in\Gamma$, let $\operatorname{tr}_{\gamma}:R(\Gamma,SO_0(4,1))\to \mathbf R$ denote the trace function.
Since  $\operatorname{tr}_{\gamma}$ is constant on orbits by conjugation, $d\operatorname{tr}_{\gamma}:Z^1(\Gamma,\mathfrak{so}(4,1))
\to \mathbf R$ vanishes on $B^1(\Gamma,\mathfrak{so}(4,1))$, and it induces a linear map,
 $d\operatorname{tr}_{\gamma}:H^1(\Gamma,\mathfrak{so}(4,1))
\to \mathbf R$. (See also~\cite{ScannellBart} for properties of the trace function)

\begin{lem}
\label{lem:R31inker} 
For a representation $\rho\! :\!\Gamma\to SO_0(3,1)\subset SO_0(4,1)$, it holds
  $$H^ 1(\Gamma,\mathbf R^{3,1})\subseteq \ker d \operatorname{tr}_{\gamma} ,\qquad  \forall \gamma\in\Gamma.$$
\end{lem}

 \begin{proof}
Using the embedding  $SO(3,1)\subset SO(4,1)$ of equation (\ref{eqn:canonicalinclusion}), the corresponding
embedding of $\mathbf R^{3,1}$ in $\mathfrak{so}(4,1)$ maps the vector $v\in\mathbf R^{3,1}$ to the matrix
$$
\begin{pmatrix}
0 & v \\
J v^t & 0
\end{pmatrix}.
$$
Then a path of representations $\rho_t$ tangent to a vector in 
$H^ 1(\Gamma,\mathbf R^{3,1})$ can be written, up to fist order and up to conjugation as
$$
\rho_t(\gamma)=
\left(\operatorname{Id}+ t
\begin{pmatrix}
0 & v \\
J v^t & 0
\end{pmatrix}
\right)
\begin{pmatrix}
\rho(\gamma) & 0 \\
0 & 1
\end{pmatrix}
+o(t^2)
=
\begin{pmatrix}
\rho(\gamma) & t v \\
-t J v^t \rho(\gamma)  & 1
\end{pmatrix}+ o(t^2),
$$
$\forall\gamma\in\Gamma$.
Hence $\operatorname{tr}(\rho_t(\gamma))=
\operatorname{tr}(\rho_0(\gamma))+ o(t^2)$ and 
$\frac{d\phantom{t}}{dt}\operatorname{tr}(\rho_t({\gamma}))\vert_{t=0}=0$.
\end{proof}

Recall that $\operatorname{Isom}^+(\mathbf H^3)\cong PSL_2(\mathbf C)$. This isomorphism comes from identifying the conformal sphere $S^2=\partial_{\infty}\mathbf H^3$ with the projective line $\mathbf {P^1C}$. 
The relation between traces in $SO(3,1)$ and $SL_2(\mathbf C)$ is the following one.

\begin{rem} Let $A\in \operatorname{Isom}^+(\mathbf H^3)$. Then:
$$
\operatorname{tr}_{SO(3,1)}(A)=\vert \operatorname{tr}_{SL_2(\mathbf C)}(A)\vert ^2
\qquad\textrm{and}\qquad 
\operatorname{tr}_{SO(4,1)}(A)=\vert \operatorname{tr}_{SL_2(\mathbf C)}(A)\vert ^2+1.
$$
\end{rem}

\begin{lem}
\label{lem:holsl2}
Let $\phi:\pi_1(M)\to PSL_2(\mathbf C)$ be the holonomy representation of a hyperbolic manifold with one cusp.
Then:
\begin{itemize}
\item[(a)] $H^1(M, sl_2(\mathbf C)_{Ad\phi})\cong T_{\phi}X(M,PSL_2(\mathbf C))\cong \mathbf C$.
\item[(b)] For every $\gamma\in\pi_1(\partial M)$, $\gamma\neq 1$,  $d \operatorname{tr}_{\gamma}:
H^1(M, sl_2(\mathbf C)_{Ad\phi})\to \mathbf C$ is nonzero.
\item[(c)] Given $\gamma_1,\gamma_2\in\pi_1(\partial M)$ so that $\phi (\gamma_i)=\pm 
\begin{pmatrix}
 1& x(\gamma_i) \\ 0 & 1
\end{pmatrix}$, then
$$
\frac{d \operatorname{tr}_{\gamma_1}}{d \operatorname{tr}_{\gamma_2}} = \left( 
\frac{x(\gamma_1)}{x(\gamma_2)}
\right)^2.
$$
\end{itemize}
\end{lem}

\begin{proof}
Assertions (a) and (b) are the contents of hyperbolic Dehn filling theorem,
cf.~\cite{KapovichBook}. To prove (c), we follow \cite[Appendix B]{BP},
and of course Thurston's notes \cite{ThurstonNotes}.  
We write a deformation as
$\phi_t(\gamma_i)=\pm \begin{pmatrix}
 e^{u_i(t)}& x(\gamma_i)(t) \\ 0 & e^{-u_i(t)}
 \end{pmatrix} $.
The parameter $t$ is not a parameter of the deformation space because the trace function
\begin{equation}
\operatorname{tr}_{\gamma_i}=  e^{u_i(t)} +  e^{-u_i(t)}=2\cosh(u_i(t))
\end{equation}
defines a local parameter of $R(M,SL_2(\mathbf C))/\! /SL_2(\mathbf C)$. However
one can take $t$ and $u_i(t)$ as analytic functions, by working in a
double branched covering, cf.~\cite{BP}. 
Following Thurston's notes, the commutativity relation becomes:
$x(\gamma_1) \sinh u_2=x(\gamma_2) \sinh u_1$.
Thus 
\begin{equation}
\lim_{t\to 0}\dot u_1/\dot u_2=\lim_{t\to 0} u_1/u_2= x(\gamma_1)(0)/x(\gamma_2)(0),
\end{equation}
and (c) follows from straightforward computation.
\end{proof}

\section{The variety of representations around $\rho_0$}
\label{sec:localanalysis}

In this section we study the geometry of $R(M,SO(4,1))$ in a neighborhood of
the holonomy representation $\rho_0$ for the complete structure of $M$, namely the 
differentiable tangent cone in Subsection~\ref{subsectio:differentialtgcone}, and
the partial slice in Subsection~\ref{subsection:partialslice}. Both tools are going 
to be used in the proof of Theorems~\ref{thm:localrigidity} and~\ref{prop:limitslope}.

\subsection{The differentiable  tangent cone}
\label{subsectio:differentialtgcone}

There are several notions of tangent cone. For our purpose, we consider the following one:

\begin{dfn}
The \emph{differentiable tangent cone} of the variety of
representations at $\rho_0$ is defined to be the subset of vectors in 
the Zariski tangent space $Z^ 1(M,\mathfrak{so}(4,1))$ that  are the tangent vector to a curve of representations, parametrized by $[0,\varepsilon)$
 and which is 
differentiable to the right at the origin.
 \end{dfn}

As the cocycles project to  cohomology, we look at the image
of the differentiable tangent cone in 
$$
H^ 1(M,\mathfrak{so}(4,1))\cong H^ 1(M,\mathfrak{so}(3,1))\oplus H^ 1(M,\mathbf R^{3,1}).
$$

The aim of this section is to prove:
\begin{prop}\label{prop:tangentcone}
If $PH^1(M,\mathfrak{so}(4,1))=0$, then the image of the differentiable tangent cone in cohomology is
contained in 
$H^ 1(M,\mathfrak{so}(3,1))\cup
 H^ 1(M,\mathbf R^{3,1})$.
\end{prop}

Before, we need a couple of remarks
and a lemma. The remarks
 can be proved by straightforward computation.

\begin{rem}
\label{rem:traceSO(4,1)}
 Let  $A\in SO_0(4,1)$. 
\begin{itemize}
 \item[--] If $A$ is parabolic and induces a translation in $\partial \mathbf H^4\setminus \{\infty\}\cong \mathbf R^3$,
 then $\operatorname{trace}(A)=5$.
 \item[--] If $A$ is parabolic and induces a screw motion  of angle $\alpha$ in $\partial \mathbf H^4\setminus \{\infty\}$,
 then $\operatorname{trace}(A)=3+2\cos(\alpha)$.
 \item[--] If $A$ is elliptic with rotational angles $\alpha$ and $\beta$, then 
$$
\operatorname{trace}(A)=1+2\cos(\alpha)+2\cos(\beta).
$$
 \item[--] If $A$ is loxodromic with translation length $\lambda$ and with angle $\alpha$, then
$$\operatorname{trace}(A)=1+2\cosh(\lambda)+2\cos(\alpha).$$
 \end{itemize}
In particular, if $A\in SO_0(4,1)$ satisfies $\operatorname{trace}(A) > 5$, then $A$ is loxodromic.
 \end{rem}

The following lemma can be found in \cite[Lemma 4.2]{ScannellSurgery}.

\begin{lem}
\label{cor:loxodromicfuchsian}
A representation of $\mathbf Z\times\mathbf Z$ in $SO_0(4,1)$ containing loxodromic elements
and obtained by perturbing a parabolic one (that consists only of translations in $\mathbf R^3$) must be conjugated to  $SO_0(3,1)$.
\end{lem}

The next lemma is the analogue of Proposition~\ref{prop:tangentcone}
for $\pi_1(\partial M)\cong\mathbf Z\oplus\mathbf Z$, 
instead of $M$.

\begin{lem}\label{cor:exclaim}
  Let $c:[0,1]\to R(\pi(\partial M),SO(4,1))$ be a
path differentiable to the right at $0$, and 
 such that $c(0)$ is the restriction of the
  holonomy of the complete structure on $M$. 
Assume that the $\mathfrak{so}(3,1)$-component of the projection 
of $\dot c(0)$ in 
$H^1(\partial M,\mathfrak{so}(3,1))\subset H^1(\partial M,\mathfrak{so}(4,1))$ is not zero.

Then,
there is $\varepsilon >0$ such that for all $0<t<\varepsilon$, $c(t)$ is
conjugate to $SO(3,1)$. In particular, the $\mathbf
R^{3,1}$-component of $\dot c(0)$ vanishes.
\end{lem}

\begin{proof} Let us identify $\dot c(0)$ with its projection in cohomology. 
Let $u$ be the $\mathfrak{so}(3,1)$ component of $\dot c(0)$, 
i.e.\ $\dot c(0)-u\in H^1(\partial M, \mathbf R^{3,1})$.  
By Lemma~\ref{lem:R31inker}, 
$$
\forall\gamma\in\pi_1(\partial M),\qquad d \operatorname{tr}_{\gamma}(\dot c(0))=d \operatorname{tr}_{\gamma}(u).
$$

As a representation in $PSL_2(\mathbf C)$, $c(0)$ is conjugated to
$$
\gamma\mapsto \pm\begin{pmatrix} 	1 & x(\gamma) \\ 
					0 & 1
                 \end{pmatrix}, \qquad
\forall \gamma\in\pi_1(\partial M).
$$
The map $x:\pi_1(\partial M)\to \mathbf C$ is a morphism and
it defines a lattice of $\mathbf C$.  
From the formula
$$
 \operatorname{tr}_{SO(4,1)}= \operatorname{tr}_{SO(3,1)}+1
= \vert \operatorname{tr}_{SL(2,\mathbf C)}\vert ^2+1
$$
and taking a lift such that $\operatorname{tr}_{SL(2,\mathbf C),\gamma}(c(0))=2$, we get
$$
d \operatorname{tr}_{SO(4,1),\gamma}(u)= 4 \operatorname{Re} (d \operatorname{tr}_{SL(2,\mathbf C),\gamma}(u)).
$$
Since  $x:\pi_1(\partial M)\to \mathbf C$ is a lattice, 
 the set of arguments 
$$
\left\{\frac{x(\gamma)}{\vert x(\gamma)\vert}\in S^1\mid 
\gamma\in\pi_1(\partial M)\setminus\{1\}\right\}
$$
is dense in the unit circle $S^1$.
By Lemma~\ref{lem:holsl2} (c), and since $ d
\operatorname{tr}_{SL(2,\mathbf C),\gamma}(u)\neq 0$, 
one can find $\gamma$ such that $\operatorname{Re} (d
\operatorname{tr}_{SL(2,\mathbf C),\gamma}(u))>0$. 

Therefore, for small enough $t$ the representation $c(t)$ contains
loxodromic elements because of Remark~\ref{rem:traceSO(4,1)},
 and is conjugated to $SO(3,1)$ by 
Lemma~\ref{cor:loxodromicfuchsian}.
\end{proof}

\begin{proof}[Proof of Propostion~\ref{prop:tangentcone}]
Since $PH^1(M,\mathfrak{so}(4,1))$ vanishes, 
the restriction from $\pi_1(M)$ to $\pi_1(\partial M)$ induces an inclusion in cohomology
$$
0\to H^1( M,\mathfrak{so}(4,1))\to H^1(\partial M,\mathfrak{so}(4,1)),
$$ and it is sufficient to 
show the corresponding statement for $\partial M$ instead of $M$.
Namely, we have  to check that ``mixed'' elements $(u,v)\in H^
1(\partial M,\mathfrak{so}(3,1) ) 
\oplus H^1(\partial M,\mathbf R^{3,1 })$ 
with $u\neq 0$ and
$v\neq 0$ are not contained in the image of the differentiable tangent
cone, and this is true by Lemma~\ref{cor:exclaim}.
\end{proof}

\begin{rem}
Further work would yield that the inverse image of the subspace $H^1(M,\mathbf R^{3,1})\subset H^1(M,\mathfrak{so}(4,1))$
in $ Z^1(M,\mathfrak{so}(4,1))$ is integrable, using the computations
of Subsection~\ref{subsection:rotations} and an argument analogue to 
\cite[Thm.~9.4]{BLP}.
\end{rem}

\subsection{The partial slice}
\label{subsection:partialslice}

Instead of working with the space of conjugacy classes of representations in $R(M,SO(4,1))$,  we shall construct
a partial slice to the orbit by conjugation. Since $SO(4,1)$ is a real group, the space of conjugacy classes 
is not an algebraic variety, though the partial slice is.

For $\rho\in R(M,SO(4,1))$, the orbit of $\rho$ by conjugation is
denoted by $O(\rho)$.

\begin{prop}
\label{prop:slice}
There exists an algebraic subvariety $S\subset R(M, SO(4,1))$ and  an
open neighborhood $U\subset R(M, SO(4,1))$ of $\rho_0$, with the
following properties: 
\begin{enumerate}
 \item $S\cap O(\rho_0)=\{\rho_0\}$. 
  \item If $\rho\in U$ satisfies that $\rho\vert_{\partial M}$ fixes a point
in $\partial \mathbf H^4$, then $O(\rho)\cap S\cap U \neq\emptyset$ and
consists of a single point. 
\item The map $TS\to H^1(M,\mathfrak{so}(4,1))$ is injective. 
\item If  $\rho\in U$ satisfies that  $O(\rho)\cap R(M,SO(3,1)) \neq \emptyset$,
then $O(\rho)\cap S\cap U\subset R(M,SO(3,1))$.
\end{enumerate}
\end{prop}

\begin{proof}
Let $\gamma_1,\gamma_2$ be a pair of peripheral elements that generate a non elementary subgroup
of $\pi_1(M)$.
Define $S$ as the subset of representations $\rho$ in $R(M,SO(4,1))$ that satisfy the conditions (a)-(d) below:
\begin{itemize}
 \item[(a)] $\rho(\gamma_1)$ fixes the same point of $\partial\mathbf H^4$ as $\rho_0(\gamma_1)$.
 \item[(b)] $\rho(\gamma_2)$ fixes the same point of $\partial\mathbf H^4$ as $\rho_0(\gamma_2)$.
\end{itemize}
Let $\mu_1\in\pi_1(M)$ be an element such that $\gamma_1$ and $\mu_1$ generate a peripheral subgroup $\pi_1(\partial M)$.
We fix an identification of $\partial \mathbf H^4\setminus Fix(\rho_0(\gamma_1))$ with $\mathbf R^3$,
so that $\rho\in S$ restricted to $\mathbf R^3$ acts by affine transformations.
Let $\vec 0\in\mathbf R^3$ denote the origin, $\rho_0(\gamma_1)(\vec 0)\neq\vec 0$.
The remaining conditions defining $S$ are:
\begin{itemize}
 \item[(c)]  $\rho(\gamma_1)(\vec 0)=\rho_0(\gamma_1)(\vec 0)$.
\item[(d)] The ordered pairs of vectors $(\rho(\gamma_1)(\vec 0),\rho(\mu_1)(\vec 0))$ and
$(\rho_0(\gamma_1)(\vec 0),\rho_0(\mu_1)(\vec 0))$ span the same oriented plane of $\mathbf R^3$.
\end{itemize}

Given a representation such that its restriction to ${\partial M}$ fixes a point in $\partial \mathbf H^4$, conditions (a) to (d) can be achieved by conjugation. Notice also that this determines the representation up to conjugacy. Namely,
conditions (a) and (b) fix a representation of the conjugacy class, up to  isometries that preserve a pair of points in
$\partial\mathbf H^4=\mathbf R^3\cup \{\infty\}$. Hence we 
may assume that the fixed points are $\{\infty\}$ for  $\rho(\gamma_1)$ in (a) and
$\vec 0$ for  $\rho(\gamma_2)$ in (b). Thus the group of elements that fix those points is the product of the orthogonal group with the group of homotethies in $\mathbf R^3$, but this indeterminacy is
eliminated by (c) and (d). hence assertions 1 and 2 of the proposition follow.

The restrictions (a) - (d) can be written as $F^{-1}(0)$ for some map $F: U\to \mathbf R^{10}$ transverse to the orbit $O(\rho_0)$, hence
$O(\rho_0)$ and $S$ are transverse. Since $B^1(M,\mathfrak{so}(4,1))$ is the tangent space to $O(\rho_0)$, $B^1(M,\mathfrak{so}(4,1))\cap T_{\rho_0}S=0$, and assertion 3 follows.
Finally, assertion 4 holds because uniqueness of assertion 2, and the fact that properties  (a) - (d) may be achieved by conjugation in $SO(3,1)$.
\end{proof}

Notice that representations such that its restriction to $\partial M$ is contained in $SO(4)$ are excluded by this set $S$, this is why we call it \emph{partial} slice.

\section{Non-Fuchsian representations}
\label{sec:sequencesreps}

This last section is devoted to the proof of Theorem~\ref{prop:limitslope}.
In Subsection~\ref{subsec:defmsscrew} we deform a parabolic group of translations in the plane
as a group of screw motions in Euclidean space. 
Viewing $\mathbf R^3$ as $\partial\mathbf H^4\setminus\{\infty\}$, those give infinitesimal
 deformations in
$\mathfrak{so}(4,1)$ that take values in $\mathbf R^{3,1}$. We claim in Lemma~\ref{lem:rotofeverycocycle}
that those are all possible deformations of $\mathbf Z\oplus \mathbf Z$ that take values in 
$\mathbf R^{3,1}$. Subsection~\ref{subsect:sequencesnonfuchsian} is devoted to the proof of 
Theorem~\ref{prop:limitslope}, assuming Lemma~\ref{lem:rotofeverycocycle}   which is proved in
Subsection~\ref{subsection:rotations}.

\subsection{Deformations with peripheral screw motions}
\label{subsec:defmsscrew}

In this section we construct explicit examples of deformations of a
parabolic representation of $\mathbf Z\oplus\mathbf Z$ that give 
cocycles valued in $\mathbf R^{3,1}$. The aim will be to show later that those are
all  the possible infinitesimal non-Fuchsian deformations.

\begin{ex} \emph{Translation along a line as limits of rotations in the plane}.
\label{ex:line}
Consider the translation that maps $x\in\mathbf R$ to $x+a$.
We extend it to a translation of the plane with vector $\begin{pmatrix}
                                                         a \\ 0
                                                        \end{pmatrix}
\in\mathbf R^2$.
For $a\in\mathbf R$, consider the family of rotations of $\mathbf R^2$ 
parametrized by $0<t<\varepsilon$,
centered at the point $\begin{pmatrix} 0 \\ 1/t \end{pmatrix} $ and of angle $\alpha= a\, t$. They can be written as:
\begin{equation}
 \begin{pmatrix}
 x \\ y
\end{pmatrix}
\mapsto 
\begin{pmatrix}
 \cos(a\, t) & - \sin (a\, t) \\
 \sin(a\, t) &  \cos(a\, t)
\end{pmatrix}
\begin{pmatrix}
 x \\ y
\end{pmatrix}
+ 
\begin{pmatrix}
 { \sin(a\, t)}/t \\ { (1-\cos(a\, t))}/t
\end{pmatrix},
\qquad
\forall \begin{pmatrix}
 x \\ y
\end{pmatrix}\in\mathbf R^2.
\end{equation}
Obviously, when $t\to 0$, this converges to the translation of vector 
$\begin{pmatrix}
a \\
0
\end{pmatrix}
$.
\end{ex}

We now want to compute the derivative of this expression with respect to $t$. 
Consider a family of representations of $\mathbf Z$ that map $1$ to the previous example.
The corresponding cocycle maps $1$ to the infinitesimal rotation 
$$
 \begin{pmatrix}
  0 &  - a \\  a & 0
 \end{pmatrix}
$$ 
plus a vertical infinitesimal translation. 
The corresponding Killing vector field is perpendicular to 
the horizontal coordinate axis $\{(x,0)\mid x\in\mathbf R\}$.

\begin{ex} \emph{Translation along a plane as limits of rotations in the space.}
\label{ex:plane}
 The same picture as above can be adapted for a group of translations of the 
plane by decomposing it as the orthogonal sum of two lines:
$$
\mathbf R^2=\langle\omega\rangle \oplus \langle\mathbf i\, \omega\rangle ,
$$
where $\omega\in \mathbf R^2$ satisfies $\vert\omega\vert=1$ and 
identifying $\mathbf R^2\cong \mathbf C$, $\mathbf i\, \omega$ denotes the result of rotating 
$\omega$ by $\pi/2$. Then in the direction $\omega$ we do not make any deformation, and in 
the direction of $\mathbf i\, \omega$
we do the construction of Example~\ref{ex:line}.
\end{ex}

Assume we have a representation $\phi\!:\Gamma\to \mathbf R^2$ into the group of translations 
Let $\operatorname{rot}\! : \mathfrak{Isom}(\mathbf R^3)\to \mathfrak{so}(3)$ denote the projection induced by taking 
the linear part of an isometry. There is a natural identification $\mathfrak{so}(3)\cong \mathbf R^3$.
An elementary computation then shows:

\begin{lem}
 If $d$ is the cocycle of deformation in Example~\ref{ex:plane}, then
$$
\operatorname{rot}(d(\gamma))= (\phi(\gamma)\cdot \mathbf i\, \omega ) \, \omega  \qquad\forall\gamma\in \Gamma.
$$
\end{lem}

Here the dot denotes the Euclidean scalar product, so that $\phi(\gamma)\cdot \mathbf i\, \omega$ denotes the orthogonal
projection of $\phi(\gamma)$ in the direction perpendicular to $\omega$.

We are interested  in the restriction of the holonomy $\rho_0$ of the complete structure, that we view as a representation
into the group of translations
 $\rho_0:\ \pi_1(\partial M)\to \mathbf R^2$.

\begin{lem}
\label{lem:rotofeverycocycle}
\begin{itemize}
 \item[(a)]
 Any cocycle $d \in Z^1(\partial M, \mathbf R^{3,1})$ takes values in $\mathbf R^{3,1}\cap \mathfrak{Isom}(\mathbf R^3)$.

\item[(b)] Every cohomology class is represented by a cocycle $d$, so that there exists a unit vector $\omega\in\mathbf R^2$,
and a parameter $\lambda\in\mathbf R$ such that
$$
\operatorname{rot}(d(\gamma))= (\phi(\gamma)\cdot \mathbf i\, \omega )\lambda\,  \omega  \qquad\forall\gamma\in \pi_1(\partial M).
$$ 
Moreover the cohomology class is trivial iff $\operatorname{rot}\circ d$ is trivial.
\end{itemize}
\end{lem}

The proof is postponed to Subsection~\ref{subsection:rotations}.

\subsection{Sequences of non-Fuchsian representations}
\label{subsect:sequencesnonfuchsian}

For all $n\in\mathbf N$, let $\rho_n:\pi_1(M)\to SO(4,1)$ 
be a representation induced by a discrete and faithful representation of 
$M_{p_n/q_n}$ not conjugated to $SO(3,1)$.

We claim that the restriction of $\rho_n$ to $\pi_1(\partial M)$ fixes
a point in $\partial \mathbf H^4$, because its restriction to
$\pi_1(\partial M)$ is a discrete and faithful representation of
$\mathbf Z$, and therefore it cannot fix an interior point 
of $\mathbf H^4$. Thus we can apply Proposition~\ref{prop:slice}, and
assume, possibly up to conjugation, that $\rho_n$ belongs to the partial
 slice $S$.

From now on, all statements are up  to subsequence. This
can be done because the limit $l$ will depend only on $M$.

\begin{lem}
\label{lem:tgvector}
There exists a cocycle $d\in T_{\rho_0}S$ and a sequence of positive real numbers $\varepsilon_n\to 0$ satisfying:
$$
\rho_n(\gamma)=(1+\varepsilon_n d(\gamma)+ o(\varepsilon_n))\rho_0(\gamma), \qquad \forall\gamma\in \pi_1M.
$$
\end{lem}

Here $o(\varepsilon_n)$ denotes a term such that $o(\varepsilon_n)/\varepsilon_n\to 0$.

\begin{proof}
This is a compactness argument. Embed $R(M,SO(4,1))$ in some $\mathbf R^N$ as an algebraic subvariety, let $\varepsilon_n$
 be
the distance between $\rho_n$ and $\rho_0$ and take a converging subsequence of unitary vectors
$\frac1{\varepsilon_n}(\rho_n-\rho_0)$. The limit must be a vector Zariski tangent to $S$, and therefore 
it is a cocycle $d$.
\end{proof}

\begin{lem}\label{lemma:5.2}
The cocycle $d$ projects to a nontrivial element in $H^1(M,\mathfrak{so}(4,1))$ contained in  
$H^1(M,\mathbf R^{3,1})$.
\end{lem}

\begin{proof} 

By Proposition \ref{prop:slice} (Assertion 3) the cocycle $d$ is non
trivial in cohomology. 
Assume that it is not contained in $H^1(M,\mathbf R^{3,1})$, we look for
a contradiction by applying the curve selection lemma for
semialgebraic sets. 

Let $S\subset R(M,SO(4,1))$ denote the slice of
Proposition~\ref{prop:slice}. 
Working with an embedding of $S\subset R(M,SO(4,1))$ in the Euclidean
space $\mathbf R^N$ and putting $\rho_0$ as the origin, let $\alpha>0$
denote the angle between $d$ 
and the linear subspace $Z^1(M, \mathbf R^{3,1})$. 
Consider the
semialgebraic cone $C$ consisting of those vectors of $\mathbf R^N$ 
whose angle with $Z^1(M, \mathbf R^{3,1})$ is $\geq \alpha/2$.
By the curve selection lemma 
applied to $C\cap S\setminus R(M,SO(3,1))$,
there exists a semialgebraic curve
$
c:[0,1]\to  S
$
such that:
\begin{enumerate}
\item $c(0)=\rho_0$.
\item $c((0,1])\subset C\cap S \setminus R(M,SO(3,1))$.
\end{enumerate}
The first non-trivial derivative $c^{(n)}(0)$ gives
an element of the
differential tangent cone whose projection to $H^1(M, \mathfrak{so}(3,1))$ is
nontrivial, by the choice of $\alpha$. Thus,
by Proposition~\ref{prop:tangentcone} applied to $c(t^{1/n})$, the cohomology class of $d$ must be 
contained in $H^1(M, \mathfrak{so}(3,1))$.

Now, we argue with the inclusion of $\partial M$ in $M$ and the
projection of $SO(4,1)$ to the de Sitter space $\mathbf
S^{3,1}=SO(4,1)/SO(3,1).$

On one hand, by Lemma~\ref{cor:exclaim} the restriction of $c$ to
$\partial M$ gives a path  
$c_\partial:[0,1]\to R(\partial M,SO(4,1))$ that must be
contained in $SO(3,1)$. In particular, the projection $\bar
c_\partial$ of such path to
$\mathbf S^{3,1}$ is the trivial path. 

On the other hand, since $c$ is not contained in $R(M,SO(3,1))$, its
projection $\bar c$ to $\mathbf S^{3,1}$ must have some non-trivial
derivative. This defines, by Corollary~\ref{cor:cyclecbar}, 
a non-trivial cocycle $b\in Z^1(M,\mathbf R^{3,1})$, which is non-trivial in
$H^1(M,\mathbf R^{3,1})$ because in the proof of Proposition~\ref{prop:slice},
it is shown that $T_{\rho_0}S\cap B^1(M,\mathfrak{so}(4,1))=0$, i.e. the tangent space to the partial
slice $S$ and the coboundary space are transverse. 
By Lemma~\ref{lem:kerboundary}, $b$ would give a non-trivial element
when restricted to the boundary, contradicting the fact that $\bar
c_\partial$ is the constant path.
\end{proof}

\begin{proof}[{Proof of Theorem~\ref{thm:localrigidity}}]
By contradiction, assume that there exists an infinite sequence ${p_n/q_n}$ so that $M_{p_n/q_n}$ 
is not locally rigid. Thus  the holonomy of $M_{p_n/q_n}$ can be perturbed 
to $\tilde\rho_n\! :\! \pi_1(M_{p_n/q_n})
\to SO(4,1)$, not contained in $SO(3,1)$. 
Notice that since the holonomy of $M_{p_n/q_n}$ 
maps the core of the filling torus to
a loxodromic element, we can assume that the restriction $\tilde\rho_n\vert_{\partial M}$ is not elliptic. Thus  $\tilde\rho_n$ belongs to the partial slice
$S$ and the arguments of  Lemmas~\ref{lem:tgvector}
 and \ref{lemma:5.2} apply. 
In other words, there is a cocycle $d$ wich is tangent to $\tilde\rho_n$, 
and that cocycle is contained in $H^1(M,\mathbf R^{3,1})$.

However, since the $\tilde\rho_n$ are obtained by perturbation of the 
fuchsian holonomies of the manifolds $M_{p_n/q_n}$, and since the sequence 
of such holonomies  
defines a cocycle in $H^1(M,\mathfrak{so}(3,1))$, we can choose the 
perturbations small enough so that the cocycle $d$ is not 
contained in $H^1(M,\mathbf R^{3,1})$. 
Hence we get a contradiction.
\end{proof}

\begin{proof}[{Proof of Theorem~\ref{prop:limitslope}}]
Let $\rho_n:\pi_1(M)\to SO(4,1)$ be a sequence of representations induced by 
a discrete and faithful representation of $\pi_1(M_ {p_n/q_n})$,
and let and $d$  be as in Lemma~\ref{lem:tgvector}.
 Since $d$ is tangent to $S$, we may assume that its restriction
to $\pi_1(\partial M)$  
takes values in the infinitesimal isometries of the Euclidean space
 $\mathbf R^3=\partial\mathbf H^4\setminus \{p_0\}$ by Lemma~\ref{lem:rotofeverycocycle} (a). Thus, the composition with the projection $\mathfrak{Isom}(\mathbf
 R^3)\to \mathfrak{so}(3)$ gives a cocycle valued in 
 infinitesimal rotations $\delta\! :\! \pi_1(\partial M)\to \mathfrak{so}(3)$.

Since $d$ is not trivial in $H^1(\partial M, \mathbf R^{3,1})$, 
by Lemma~\ref{lem:rotofeverycocycle} (b),
the image
of $\delta$ is 
non trivial and it has
an invariant direction 
$\omega \in\mathbf R^2\subset \mathbf R^3$, $\vert \omega \vert=1$.

The restriction of $\rho_n$ to $\partial M$ consists of screw motions
of $\mathbf R^3$, and for each element $\gamma\in\pi_1(\partial M)$,
the  translation length of this screw motion is the product:
$$
\operatorname{trans}(\rho_n(\gamma))=(\rho_n(\gamma)(0)-(0))\cdot \omega_n
$$
where $\omega_n\in \mathbf R^3$ is a unitary vector in the direction of the axis of the screw motion.
Since the projection of $d$ is $\delta$, 
it follows that  $\omega_n\to\omega$.  Thus,
 given a system of peripheral generators $\langle\gamma_1,\gamma_2\rangle$, and the restriction 
$$
p_n\, \operatorname{trans}(\rho_n(\gamma_1)) +q_n\, \operatorname{trans}(\rho_n(\gamma_2))=0,
$$
we deduce the limit
$$
\lim_{n\to\infty}\frac{p_n}{q_n}=-\frac{(\rho_0(\gamma_2)(0)-(0))\cdot \omega}{(\rho_0(\gamma_1)(0)-(0))\cdot \omega}
$$
which is a well defined element $l\in \mathbf R\cup\infty$, depending
only on the cusped manifold $M$. 
\end{proof}

\begin{proof}[{Proof ot Theorem \ref{thm:limitslope}}]
Let $C$ be the set of coefficients $(p,q)$ so that $\pi_1(M_{p/q})$ 
has a discrete 
and faithful representation in $SO(4,1)$, other than the holonomy of 
the complete hyperbolic structure of $M_{p/q}$. By Theorem~\ref{t_c}, 
any sequence of pairwise dinstinct such representations must converge to the 
holonomy of the complete hyperbolic structure of $M$.   
By Theorem~\ref{prop:limitslope} it now follows that 
the set $C$ is asymptotic to the line 
$p/q=l$, where $l$ is a number 
--- possibly $\infty$ --- depending only on $M$. Thus $C$ 
cannot be co-finite in the set of all filling coefficients.  
\end{proof}

\subsection{Cohomology of $\mathbf Z\oplus \mathbf Z$ with coefficients in $\mathbf R^{3,1}$}
\label{subsection:rotations}

The aim of this paragraph is to prove Lemma~\ref{lem:rotofeverycocycle}.

Before the proof we fix some notation.
The restriction of the holonomy representation $\rho_0$ of $\pi_1(M)$ to 
$\pi_1(\partial M)\cong \mathbf Z\oplus\mathbf Z$  
is a parabolic representation. 
Identifying the fixed point $p_0$ of $\rho_0|_{\pi(\partial M)}$ with
$\infty$, so that $\partial\mathbf H^3=\mathbf R^2\cup \{ p_0\}$,
the restriction is a representation by translations, that defines a
lattice in the plane $\mathbf R^2$.

We choose
$$
p_0=\begin{pmatrix} 1 \\ 1 \\ 0 \\ 0 \end{pmatrix}
$$
to be the point of the light cone invariant by the holonomy of
$\partial M$. 
With that choice, for $\gamma\in\pi_1(\partial M)$
if the translation vector of $\gamma$ is 
$$\operatorname{trans}(\rho_0(\gamma))=(x,y)\in \mathbf R^ 2,
$$
 then the holonomy (as an element of $SO(3,1)$) is:
$$
\rho_0(\gamma)=\exp\begin{pmatrix} 0 & 0 & x & y \\  0 & 0 & x & y \\ x & -x & 0 & 0 \\  y & -y & 0 & 0 \end{pmatrix}=
 \begin{pmatrix} 1 + ({x^ 2+y^ 2})/2 & -({x^ 2+y^ 2})/2 & x & y \\  ({x^ 2+y^ 2})/2 & 1-({x^ 2+y^ 2})/2 & x & y \\ x & -x & 1 & 0 \\  y & -y & 0 & 1 \end{pmatrix}.
$$

Elements in $\mathbf R^{3,1}$ may be written as
\begin{equation}
\label{eqn:coordinatesR31}
v=\begin{pmatrix} z+\lambda \\ z-\lambda \\ -\beta \\ \alpha \end{pmatrix}
\quad\textrm{ with } \lambda, z, \beta, \alpha\in\mathbf R.
 \end{equation}
Using the inclusion $\mathbf R^{3,1}\subset \mathfrak{so}(4,1)$ of Equation~\ref{eqn:r41inso41}, the parameter $\lambda$ corresponds to the length of an infinitesimal displacement of $p_0$ in $\partial\mathbf H^4$ in the direction perpendicular to
$\partial \mathbf H^3$.

 In particular  $\lambda=0$
defines the subspace of infinitesimal isometries that fix $p_0$, and therefore they restrict to infinitesimal similarities
of $\mathbf R^3=\partial \mathbf H^4\setminus \{ p_0\}$. 
Since our deformations
vanish on the direction tangent to $\mathbf R^2\times\{0\}$ they must be
infinitesimal isometries
of Euclidean space. In fact we have:
$$
  \mathbf R^{3,1}\cap \mathfrak{Isom}(\mathbf R^3)= \{ v\in\mathbf R^{3,1}\mid \lambda=0 \}.
$$
The parameter $z$ describes the length of an infinitesimal translation   in the direction perpendicular to $\mathbf R^2$.
Finally, $\beta$ and $\alpha$ correspond to an infinitesimal rotation of vector $(\alpha,\beta,0)\in \mathbf R^2\times\{ 0 \}$.
Using the coordinates in Equation~\ref{eqn:coordinatesR31},
the projection $ \operatorname{rot}\! :\mathfrak{Isom} (\mathbf R^3)\to \mathfrak{so}(3)\cong\mathbf R^3$ 
restricts to:
$$
\begin{array}{rcl}
  \operatorname{rot}: \{\lambda=0\}\subset \mathbf R^{3,1} & \to & \mathbf R^2 \times\{ 0 \}\\
v & \mapsto & (\alpha,\beta, 0)
\end{array}.
$$
Here  $\operatorname{rot}$ denotes  the  tangent map of
the epimorphism $\operatorname{Isom}(\mathbf R^3) \to O(3)$.

\begin{proof}[Proof of Lemma~\ref{lem:rotofeverycocycle}]

Fix a system of generators $g_1,g_2$ for $\pi_1(\partial M)$, so that $\rho_0(g_1)$ is a translation of vector 
$(x_1,y_1)\in\mathbf R^2$, and $\rho_0(g_2)$, of vector  $(x_2,y_2)\in\mathbf R^2$.

For a cocycle $d\in Z^1(\partial M;\mathbf R^{3,1})$, define
$\lambda_1$, $\lambda_2$, $\alpha_1$, $\alpha_2$, $\beta_1$, $\beta_2$, $z_1$ and $z_2\in\mathbf R$
so that 
$$
d(g_i)=\begin{pmatrix}
        z_i+\lambda_i \\ z_i-\lambda_i \\ -\beta_i \\ \alpha_i
       \end{pmatrix}, \qquad \textrm{ for } i=1,2.
$$
By Fox calculus, cf.\ \cite{LubotzkyMagid}, the parameters $\alpha_i$, $\beta_i$, $z_i$ and
$\lambda_i$ are subject to the relation $(g_1-1) d(g_2)=(g_2-1) d(g_1)$, which is equivalent to
$$
\lambda_1=\lambda_2=0, \qquad -\beta_2 \, x_1+\alpha_2\, y_1=-\beta_1\, x_2 + \alpha_1\, y_2,
$$
that give a five dimensional real space on the parameters 
subject to these relations.
In particular since $\lambda_1=\lambda_2=0$, statement (a) is proved.

Notice that $d$ is a coboundary if and only if there exist parameters $A, B, L\in\mathbf R$ such that
$$
d(g_i)=\begin{pmatrix}
		L (x_i^2+y_i^2)-B\, x_i+A\, y_i  \\ 
		L (x_i^2+y_i^2)-B\, x_i+A\, y_i  \\
		L \, 2\, x_i \\
		L \, 2\, y_i
       \end{pmatrix}, \qquad \textrm{ for } i=1,2.
$$

It follows immediately that if $\alpha_1=\alpha_2=\beta_1=\beta_2=0$, then $d$ is a coboundary.

Now the remaining of the proof is an elementary but tricky computation. The equality
$$
-\beta_2 \, x_1+\alpha_2\, y_1=-\beta_1\, x_2 + \alpha_1\, y_2,
$$
may be seen as an inequality of imaginary parts:
$$
\Im ( (x_1-\mathbf i\, y_1) (\alpha_2+\mathbf i\, \beta_2)) =
\Im ( (x_2-\mathbf i\, y_2) (\alpha_1+\mathbf i\, \beta_1)). 
$$
Next we claim that, by adding a cocycle, we can remove imaginary parts.
Namely the expression 
\begin{equation}
 \label{eqn:tovanish}
 (x_1-\mathbf i\, y_1) (\alpha_2+\mathbf i\, \beta_2) -
  (x_2-\mathbf i\, y_2) (\alpha_1+\mathbf i\, \beta_1) 
\end{equation}
may have nontrivial real part, but we can assume that it vanishes,
because adding the cocycle such that
$\alpha_j = L y_j$ and $\beta_j= - L x_j$ for some $L\in\mathbf R$ and $j=1,2$, it means 
changing 
the expression (\ref{eqn:tovanish})
by adding  
$2\, L (x_1y_2-y_1x_2)\neq 0$.

Since (\ref{eqn:tovanish}) vanishes, there exist $\lambda\in\mathbf R$ and $\omega\in\mathbf C$
with $\vert\omega\vert=1$ such that :
$$
\frac{\alpha_2+\mathbf i\, \beta_2}{x_2-\mathbf i\, y_2}=
\frac{\alpha_1+\mathbf i\, \beta_1}{x_1-\mathbf i\, y_1}=
\frac\lambda{2}\,\mathbf i\, \omega^2.
$$
Adding the coboundary such that $\alpha_j = \frac\lambda{2} y_j$ and $\beta_j= - \frac\lambda{2} x_j$,
for $j=1,2$,
we deduce:
$$
\begin{array}{rcl}
\alpha_1+\mathbf i\, \beta_1 & = & (x_1-\mathbf i\, y_1) \frac\lambda{2}\,\mathbf i\, \omega^2 +
              \frac\lambda{2} (y_1-\mathbf i \, x_1) \\
\alpha_2+\mathbf i\, \beta_2 & = & (x_2-\mathbf i\, y_2) \frac\lambda{2}\,\mathbf i\, \omega^2 +
              \frac\lambda{2} (y_2-\mathbf i \, x_2).
\end{array}
$$
Hence, expressing the scalar product $\cdot$ in terms of real parts, we have:
$$
\begin{array}{rcl}
(\alpha_1, \beta_1 )=\alpha_1+\mathbf i\, \beta_1 & = & \Re( (x_1-\mathbf i\, y_1) \,\mathbf i\, \omega) \lambda\, \omega =
       ((x_1,y_1)\cdot \mathbf i\,\omega)\,  \lambda\, \omega
             \\
(\alpha_2, \beta_2 )=\alpha_2+\mathbf i\, \beta_2 & = & \Re( (x_2-\mathbf i\, y_2) \,\mathbf i\, \omega)  \lambda\, \omega =
       ((x_2,y_2)\cdot \mathbf i\,\omega) \,  \lambda\, \omega
\end{array}
$$
And we conclude the proof of the lemma by linearity.
\end{proof}

\bibliographystyle{plain}

\bibliography{refsfuchsian}

\begin{thebibliography}{10}

\bibitem{ScannellBart}
Anneke Bart and Kevin~P. Scannell.
\newblock A note on stamping.
\newblock {\em Geom. Dedicata}, 126:283--291, 2007.

\bibitem{Bestvina}
Mladen Bestvina.
\newblock Degenerations of the hyperbolic space.
\newblock {\em Duke Math. J.}, 56(1):143--161, 1988.

\bibitem{BLP}
Michel Boileau, Bernhard Leeb, and Joan Porti.
\newblock Geometrization of 3-dimensional orbifolds.
\newblock {\em Ann. of Math. (2)}, 162(1):195--290, 2005.

\bibitem{BP}
Michel Boileau and Joan Porti.
\newblock Geometrization of 3-orbifolds of cyclic type.
\newblock {\em Ast\'erisque}, (272):208, 2001.
\newblock Appendix A by Michael Heusener and Porti.

\bibitem{EpsteinPenner}
D.~B.~A. Epstein and R.~C. Penner.
\newblock Euclidean decompositions of noncompact hyperbolic manifolds.
\newblock {\em J. Differential Geom.}, 27(1):67--80, 1988.

\bibitem{JohnsonMillson}
Dennis Johnson and John~J. Millson.
\newblock Deformation spaces associated to compact hyperbolic manifolds.
\newblock In {\em Discrete groups in geometry and analysis (New Haven, Conn.,
  1984)}, volume~67 of {\em Progr. Math.}, pages 48--106. Birkh\"auser Boston,
  Boston, MA, 1987.

\bibitem{KapovichArticle}
Michael Kapovich.
\newblock Deformations of representations of discrete subgroups of {${\rm
  SO}(3,1)$}.
\newblock {\em Math. Ann.}, 299(2):341--354, 1994.

\bibitem{KapovichBook}
Michael Kapovich.
\newblock {\em Hyperbolic manifolds and discrete groups}, volume 183 of {\em
  Progress in Mathematics}.
\newblock Birkh\"auser Boston Inc., Boston, MA, 2001.

\bibitem{KapovichPreprint07082671}
Michael Kapovich.
\newblock On sequences of finitely generated discrete groups.
\newblock {\em Preprint, arXiv:0708.2671}, 2007.

\bibitem{LubotzkyMagid}
Alexander Lubotzky and Andy~R. Magid.
\newblock Varieties of representations of finitely generated groups.
\newblock {\em Mem. Amer. Math. Soc.}, 58(336):xi+117, 1985.

\bibitem{Morgan}
John~W. Morgan.
\newblock Group actions on trees and the compactification of the space of
  classes of {${\rm SO}(n,1)$}-representations.
\newblock {\em Topology}, 25(1):1--33, 1986.

\bibitem{MorganShalen}
John~W. Morgan and Peter~B. Shalen.
\newblock Valuations, trees, and degenerations of hyperbolic structures. {I}.
\newblock {\em Ann. of Math. (2)}, 120(3):401--476, 1984.

\bibitem{MorganShalenIII}
John~W. Morgan and Peter~B. Shalen.
\newblock Degenerations of hyperbolic structures. {III}. {A}ctions of
  {$3$}-manifold groups on trees and {T}hurston's compactness theorem.
\newblock {\em Ann. of Math. (2)}, 127(3):457--519, 1988.

\bibitem{Paulin}
Fr{\'e}d{\'e}ric Paulin.
\newblock Sur les automorphismes ext\'erieurs des groupes hyperboliques.
\newblock {\em Ann. Sci. \'Ecole Norm. Sup. (4)}, 30(2):147--167, 1997.

\bibitem{Scannell}
Kevin~P. Scannell.
\newblock Infinitesimal deformations of some {${\rm SO}(3,1)$} lattices.
\newblock {\em Pacific J. Math.}, 194(2):455--464, 2000.

\bibitem{ScannellSurgery}
Kevin~P. Scannell.
\newblock Local rigidity of hyperbolic 3-manifolds after {D}ehn surgery.
\newblock {\em Duke Math. J.}, 114(1):1--14, 2002.

\bibitem{ThurstonNotes}
William~P. Thurston.
\newblock {The Geometry and Topology of Three-Manifolds}.

\bibitem{Weil}
Andr{\'e} Weil.
\newblock Remarks on the cohomology of groups.
\newblock {\em Ann. of Math. (2)}, 80:149--157, 1964.

\bibitem{Wolf}
Joseph~A. Wolf.
\newblock {\em Spaces of constant curvature}.
\newblock Publish or Perish Inc., Houston, TX, fifth edition, 1984.

\end{thebibliography}

\textsc{Dipartimento di Matematica Applicata ``U. Dini'', Univestit\`a
  degli studi di Pisa, via Buonarroti 1c, 56127, Italy}

\textsc{Departament de Matem\`atiques, Universitat Aut\`onoma de Barcelona, E-08193 Bellaterra, Spain}

\end{document}